\documentclass[10pt,a4paper,reqno]{amsart}

\usepackage[latin1]{inputenc}
\usepackage{amsmath}
\usepackage{amssymb,amsthm}
\usepackage[mathscr]{eucal}

\newtheorem{thm}{Theorem}[section]
\newtheorem{prop}[thm]{Proposition}
\newtheorem{lemma}[thm]{Lemma}
\newtheorem{cor}[thm]{Corollary}
\theoremstyle{definition}

\newtheorem*{remark}{Remark}

\newtheorem{defn}[thm]{Definition}
\newtheorem*{prf}{Proof}

\numberwithin{equation}{section}

\DeclareMathOperator{\supp}{supp}
\DeclareMathOperator{\ind}{ind}

\let\epsilon=\varepsilon
\let\phi=\varphi
\def\btau{\boldsymbol{\tau}}
\newcommand\sig[1]{\ensuremath{\sigma(#1)}}

\newcommand\distribD{\ensuremath{\mathcal{D}}}

\newcommand\defD{\ensuremath{\mathscr{D}}}
\newcommand\testspace{\ensuremath{C_0^\infty(\mathbb{R}^2\setminus\{0\})}}

\newcommand\dz{\ensuremath{\frac{\partial}{\partial z}}}
\newcommand\dzbar{\ensuremath{\frac{\partial}{\partial \bar{z}}}}
\newcommand\dr{\ensuremath{\frac{\partial}{\partial r}}}
\newcommand\dtheta{\ensuremath{\frac{\partial}{\partial \theta}}}

\newcommand\psip{\ensuremath{\psi_+}}
\newcommand\psim{\ensuremath{\psi_-}}

\newcommand\mpauli{\ensuremath{\mathfrak{P}_{\text{max}}}}
\newcommand\evpauli{\ensuremath{\mathfrak{P}_{\text{EV}}}}

\newcommand\pauli{\ensuremath{\mathfrak{P}}}

\newcommand\diffq{\ensuremath{q}}
\newcommand\qminus{\ensuremath{q_{-}}}
\newcommand\qplus{\ensuremath{q_{+}}}

\newcommand\dirac{\ensuremath{\mathfrak{D}}}
\newcommand\diraca{\ensuremath{\mathfrak{d}}}
\newcommand\diracmin{\ensuremath{\dirac_{\text{min}}}}

\newcommand\cldiracmin{\ensuremath{\overline{\dirac}_{\text{min}}}}
\newcommand\adcldiracmin{\ensuremath{\overline{\dirac}^*_{\text{min}}}}

\begin{document}

\title[Pauli operator with several AB solenoids]{On the Dirac and Pauli operators with\\several Aharonov-Bohm solenoids}
\author[M. Persson]{Mikael Persson}
\address[M. Persson]{Department of Mathematical Sciences\\
                        Chalmers University of Technology \\
                        and University of Gothenburg\\
                        SE-412 96 Gothenburg \\
                        Sweden}
\email{mickep@math.chalmers.se}

\maketitle

\begin{abstract}
\noindent We study the self-adjoint Pauli operators that can be realized as the square of a self-adjoint Dirac operator and correspond to a magnetic field consisting of a finite number of Aharonov-Bohm solenoids and a regular part, and prove an Aharonov-Casher type formula for the number of zero-modes for these operators. We also see that essentially only one of the Pauli operators are spin-flip invariant, and this operator does not have any zero-modes.
\end{abstract}

\vspace*{12pt}

\vspace*{12pt}

\section{Introduction}

Two-dimensional spin-$\frac{1}{2}$ non-relativistic quantum systems with magnetic fields are described by the Pauli operator. For regular magnetic fields the Pauli operator is usually defined as the square of the Dirac operator. However, for more singular magnetic fields, such as the delta field, an Aharonov-Bohm (AB) solenoid, generates (see~\cite{ab}), the situation is more complicated. Then there are many self-adjoint extensions of both the Dirac and the Pauli operator, originally defined on smooth functions with compact support not touching the singular points. The different extensions describe different physics, and it is not clear which extensions describe the real physical situation.

We consider the case of a magnetic field consisting of finitely many AB solenoids and a smooth field with compact support. Up to now only two different Pauli extensions have been studied for this type of magnetic field (see~\cite{ervo,per}), both defined via a quadratic form. Since the Pauli operator classically is the square of the Dirac operator it is motivated to study the self-adjoint Pauli extensions that can be obtained in this way. This is the purpose of the present paper.

A natural property to expect from the Pauli operator is that it transforms in an (anti)-unitarily way when the sign of the magnetic field is changed and the spin-up and spin-down components are switched. This property is usually called spin-flip invariance, and we want to answer the question of which Pauli operators defined in different ways satisfy it.

Another natural property to expect is the possibility to approximate our operator with operators corresponding to regular magnetic fields. In~\cite{bopu} it was studied which Pauli extensions are possible to approximate, and the condition was expressed in terms of the asymptotics of the functions in the domain at the singular points.

One of the ways to control the properties of self-adjoint extensions consists in the analysis of the dimension of the kernel of the operator. The dimension of the kernel of the Pauli operator is usually given by the Aharonov-Casher formula (see~\cite{ac}). The formula has been proved in different settings, see~\cite{cfks,gegr,mi}. Recently this formula was also proved for one of the extensions for a very singular magnetic field (containing the case with AB solenoids) in~\cite{ervo}. Another extension was introduced in~\cite{gegr}, and in~\cite{per} an Aharonov-Casher type formula was established for that extension.

The Dirac operator with strongly singular magnetic field has been studied before in~\cite{ar,arha,deSo,hiog,ta}. In~\cite{hiog} a formula for the dimension of the kernel of the Dirac operator was proved for two different asymmetric self-adjoint extensions (i.e. those with different behavior of spin-up and spin-down components), and it was demonstrated that, in fact, this dimension may differ for self-adjoint realizations, although, each of them seems to be quite natural. These extensions are closely related to the ones introduced in~\cite{ar}. In both these articles the magnetic field is the same as the one we consider (the one in~\cite{ar} does not have the regular part), with the addition of even more singular terms containing derivatives of the delta distributions. Using gauge transformations one can dispose of these derivatives.

In~\cite{ta} it was proved that one of these asymmetric Dirac extensions can be approximated by operators corresponding to more regular magnetic fields. However, this extension, as we show, lacks the property of being spin-flip invariant.

In Section~\ref{sec:dirac} we study the Dirac operator. In order to be able to treat the general case, we need first to repeat in details the description of all self-adjoint extensions corresponding to only one AB solenoid, given in~\cite{ta}. To extend this description to the case of several solenoids, we use the glue-together procedure, proposed in~\cite{agro}. This procedure enables one to glue together operators with different behavior at different points; however it is natural to restrict oneself to the case of the same behavior at all points. After that we check which extensions are spin-flip invariant and finally we prove a formula for the dimension of the kernel of the Dirac extensions.

In Section~\ref{sec:pauli} we consider the Pauli operators that are the square of some self-adjoint Dirac operator defined in Section~\ref{sec:dirac}. We show exactly which Pauli extensions are obtained in this way, in terms of the asymptotics of functions in the domain of the Pauli operator at the points where the singular AB solenoids are located. We also find an Aharonov-Casher type formula for these Pauli operators. It turns out that there are only two of them that have zero-modes. These two extensions are very asymmetric though, admitting singularities in one component only, which looks rather non-physical. All the other extensions have singularities in both the spin-up and spin-down components, and they are coupled.

It turns out that the Pauli operator studied in~\cite{ervo} is a mixture of these two asymmetric extensions, admitting different physical situations at different AB 
solenoids. This has to do with the normalization of the AB intensities. In that article the AB intensities were chosen to be normalized to the interval $[-1/2,1/2)$. If any 
of the intervals $(0,1)$ or $(-1,0)$ would have been chosen instead, then the Pauli operator would have been one of the asymmetric ones studied in this article. Both these 
asymmetric operators have the advantage that they describe the same physics at all AB solenoids. In the end of the article we present a discussion of the properties of the 
self-adjoint Pauli extensions with respect to different ways of normalization of AB intensities.

\section{The Dirac operator with singular magnetic field}
\label{sec:dirac}
The goal in this section is to describe the self-adjoint Dirac operators corresponding to a magnetic field consisting of several (but finitely many) AB solenoids together with a smooth field, and to find an Aharonov-Casher type formula for the dimension of the kernel of these self-adjoint operators. Let us introduce some notations that will be used throughout the article. We identify a point $(x_1,x_2)$ in $\mathbb{R}^2$ with the complex number $z=x_1+ix_2$, and we will often write $z$ in polar coordinates, $z=re^{i\theta}$. Sometimes it will be convenient to use the polar coordinates $r_je^{i\theta_j}$ with $z_j$ as the origin. The magnetic field will consist of a regular part $B_0\in C^1_0(\mathbb{R}^2)$ and a singular part consisting of $n$ AB solenoids located at the points $\Lambda=\{z_j\}_1^n$, so that the magnetic field $B$ has the form
\begin{equation}
\label{eq:magn}
B(z)=B_0(z)+\sum_{j=1}^n 2\pi \alpha_j \delta_{z_j}.
\end{equation}
Owing to gauge equivalence (see~\cite{ta}) we can assume that all the AB intensities $\alpha_j$ (fluxes divided by $2\pi$) belong to the interval $(0,1)$. All derivatives will be considered in the distribution space $\distribD'(\mathbb{C}\setminus\Lambda)$. We will denote by $h$ a magnetic scalar potential satisfying $\Delta h=B$. The magnetic scalar potential is uniquely defined modulo addition of a harmonic function. We will use the scalar potential
\begin{equation}
\label{eq:scalarpot}
h(z)=\frac{1}{2\pi}\int B_0(\zeta)\log|z-\zeta|d\lambda(\zeta)+\sum_{j=1}^n \alpha_j\log|z-z_j|=h_0(z)+\sum_{j=1}^n h_j(z),
\end{equation}
where $d\lambda$ is the Lebesgue measure. The \emph{actions} $\diffq^h_\pm$, which will be used to describe how the Dirac operator acts, are defined by
\begin{equation*}
\qplus^h u = -2i\tfrac{\partial}{\partial \bar{z}}\left(ue^{-h}\right)e^h
\end{equation*}
and
\begin{equation*}
\qminus^h u = -2i \tfrac{\partial}{\partial z}\left(ue^h\right)e^{-h}.
\end{equation*}
These actions $\diffq_+^h$ and $\diffq_-^h$ are usually called the spin-up and spin-down actions, respectively. The Dirac action is given by
\begin{equation*}
\diraca^h =
\begin{pmatrix}
0 & \qminus^h\\
\qplus^h & 0\\
\end{pmatrix}.\\
\end{equation*}
To be able to describe the self-adjoint Dirac operators with several AB solenoids we first study the self-adjoint extensions of the Dirac operator with one AB solenoid, originally defined on smooth functions with compact support not touching the singular point.

\subsection{The Dirac operator with one AB solenoid}
The case of one AB solenoid has been studied before (see~\cite{deSo,ta}), and we just sketch the way it was done since we need the detailed information about these extensions for our further analysis. We let the AB solenoid have intensity $\alpha=\alpha_1\in(0,1)$ and be located at the origin. We will describe all self-adjoint extensions of the Dirac operator originally defined on $(\testspace)^2$. Let us drop the superscript $h$ of $\diffq_\pm^h$ and $\diraca^h$ in this subsection. The actions $\diffq_\pm$ can in this case be written as
\begin{equation}
\label{eq:qplus}
\qplus = -2i\left(\dzbar-\frac{\alpha}{2\bar{z}}\right)=-ie^{i\theta}\left(\dr+\frac{i}{r}\dtheta-\frac{\alpha}{r}\right)
\end{equation}
and
\begin{equation}
\label{eq:qminus}
\qminus = -2i\left(\dz+\frac{\alpha}{2z}\right)=-ie^{-i\theta}\left(\dr-\frac{i}{r}\dtheta+\frac{\alpha}{r}\right).
\end{equation}
The minimal Dirac operator $\diracmin$, obviously symmetric, is defined by
\begin{equation*}
\begin{aligned}
\defD(\diracmin)&=(\testspace)^2,\\
\diracmin\psi&=\diraca\psi,\quad\text{for}\ \psi\in\defD(\diracmin).
\end{aligned}
\end{equation*}
It can be seen that $\cldiracmin$ has deficiency index $(1,1)$, and the deficiency spaces $\chi_\pm=\ker(\adcldiracmin\mp i)$ are spanned by 
\begin{equation}
\label{eq:defbas}
\xi_\pm(re^{i\theta})=
\begin{pmatrix}
K_{1-\alpha}(r)e^{-i\theta}\\
\pm K_{\alpha}(r)
\end{pmatrix}.
\end{equation}
Denote by $U$ any unitary operator from $\chi_+$ to $\chi_-$. Then $U$ takes $\xi_+$ to $e^{i\tau}\xi_-$ for some $\tau\in[0,2\pi)$. According to the theorem of Krein and von Neumann, described in~\cite{akgl}, all self-adjoint extensions of $\cldiracmin$ can be parameterized by $\tau$ as
\begin{align}
\label{eq:defdirac}
\defD(\dirac^\tau)&=\left\{\psi\in(L_2(\mathbb{R}^2))^2:\ \psi=\psi_0+\mu(\xi_++e^{i\tau}\xi_-),\psi_0\in\defD(\cldiracmin),\ \mu\in\mathbb{C}\right\}\\
\label{eq:dirac}
\dirac^\tau \psi&=\diraca\psi_0+i\mu(\xi_+-e^{i\tau}\xi_-),\quad\text{for}\ \psi\in\defD(\dirac^\tau).
\end{align}
It is also possible to describe the self-adjoint extensions by studying the asymptotic behavior of the functions in the domain at the origin. To see this, let us define the linear functionals $c_{-\alpha}^\pm$ and $c_{\alpha-1}^\pm$ on $\defD(\dirac^\tau)$ as
\begin{align}
\label{eq:boundary}
c_{-\alpha}^{\pm}(\psi_\pm) &= \lim_{r\to 0}r^\alpha\frac{1}{2\pi}\int_{0}^{2\pi}\psi_\pm(re^{i\theta}) d\theta,\ \text{and}\\
c_{\alpha-1}^{\pm}(\psi_\pm) &= \lim_{r\to 0}r^{1-\alpha}\frac{1}{2\pi}\int_{0}^{2\pi}\psi_\pm(re^{i\theta}) e^{i\theta} d\theta.
\end{align}
For $\psi=\psi_0+\mu(\xi_++e^{i\tau}\xi_-)$ in $\defD(\dirac^\tau)$, where $\psi_0\in\defD(\cldiracmin)$, applying these functionals gives no contribution from $\psi_0$ since the limit of functions in $\defD(\cldiracmin)$ tends to zero at the origin. Let us introduce the notation $\sig{\alpha}=\Gamma(\alpha)2^\alpha$. Using the asymptotics for the Bessel functions we get
\begin{equation}
\label{eq:funktional}
\begin{aligned}
c_{-\alpha}^+(\psip) &=0\\
c_{-\alpha}^-(\psim) &=\frac{\mu}{2}(1-e^{i\tau})\sig{\alpha}\\
c_{\alpha-1}^+(\psip) &= \frac{\mu}{2}(1+e^{i\tau})\sig{1-\alpha}\\
c_{\alpha-1}^-(\psim) &=0
\end{aligned}
\end{equation}
for such functions $\psi\in\defD(\dirac^\tau)$. Here $\mu$ is the same constant as in~\eqref{eq:defdirac}. An equivalent description of all self-adjoint Dirac extensions is
\begin{align}
\label{eq:defdiracalt}
\defD(\dirac^\tau)&=\Big\{\psi\in(L_2(\mathbb{R}^2))^2:\ 
\diraca\psi\in(L_2(\mathbb{R}^2))^2;\\
\nonumber&\quad\quad\frac{c_{\alpha-1}^+(\psip)}{c_{-\alpha}^-(\psim)}=i\cot(\tau/2)\frac{\sig{1-\alpha}}{\sig{\alpha}},\ \text{and}\\
\nonumber&\quad\quad c_{-\alpha}^+(\psip)=c_{\alpha-1}^-(\psim)=0\Big\},\\
\label{eq:diracalt}
\dirac^\tau \psi&=\diraca\psi,
\quad \text{for}\ \psi\in\defD(\dirac^\tau).
\end{align}

\subsection{The Dirac operator with several AB solenoids together with a regular field}

In this subsection we are going to study the Dirac operator for a magnetic field consisting of a finite number of AB solenoids together with a regular background field. We will use the same method as in~\cite{agro} to glue together the different self-adjoint Dirac operators corresponding to only one AB solenoid and the self-adjoint Dirac operator corresponding to the regular magnetic field.

Note here that we do not study all self-adjoint extensions but only the ones that are subject to the natural locality principle; which says that the asymptotics at a singular point of functions in the domain should be independent of the asymptotics at the other singular points.

We start by defining the Dirac operator with two AB solenoids together with a smooth field. The general case does not give any extra difficulties. Let the magnetic field $B$ consist of a smooth field $B_0$ with compact support and two AB solenoids located at $z_1$ and $z_2$ with intensities $\alpha_1$ and $\alpha_2$,
\begin{equation}
\label{eq:magntva}
B(z)=B_0(z)+2\pi\alpha_1\delta_{z_1}+2\pi\alpha_2\delta_{z_2}.
\end{equation}
In this case our scalar potential $h$ can be written as
\begin{equation*}
h(z)=h_0(z)+h_1(z)+h_2(z)=\frac{1}{2\pi}(\log|\cdot|*B_0)(z)+\alpha_1\log|z-z_1|+\alpha_2\log|z-z_2|.
\end{equation*}
From the previous section we have self-adjoint Dirac operators $\dirac^{\tau_1,h_1}$ and $\dirac^{\tau_2,h_2}$ corresponding to each of the AB solenoids separately. Let us drop the parameters $\tau_1$ and $\tau_2$ from the superscripts. So, for example, when we write $\dirac^{h_1}$ we mean some self-adjoint extension with one AB solenoid located at $z_1$.

Let $\phi_j\in C_0^\infty(\mathbb{R}^2)$, $j=1,2$, be equal to $1$ in a neighborhood of $z_j$ and have small support not touching a neighborhood of $z_k$, $k\neq j$ and $0\leq \phi_j\leq 1$. Let $\phi_0=1-\phi_1-\phi_2$. We denote by $E_{jk}$ the set $\supp\phi_j\cap\supp\phi_k$.

Let us introduce the multiplication operators $V^{h_j}$ as
\begin{equation*}
V^{h_j}=
2i
\begin{pmatrix}
0 & -\dz h_j\\
\dzbar h_j & 0\\
\end{pmatrix}.
\end{equation*}
Note that $V^{h_0}$ is bounded from $(L_2)^2$ to $(L_2)^2$. For $j\neq 0$ we will be sure to apply the operators $V^{h_j}$ only on functions being zero in a neighborhood of the singular points $z_j$.

\begin{defn}
\label{def:dirac}
The Dirac operator $\dirac^{h}$ corresponding to the magnetic field $B$ in~\eqref{eq:magntva} is defined as
\begin{equation*}
\defD(\dirac^{h})=\{\psi\in(L_2(\mathbb{R}^2))^2:\ \phi_j\psi\in\defD(\dirac^{h_j}),\ j=0,1,2\}
\end{equation*}
and
\begin{equation*}
\begin{aligned}
\dirac^{h}\psi=& (\dirac^{h_0}+V^{h_1}+V^{h_2})(\phi_0\psi)\\
+&(\dirac^{h_1}+V^{h_0}+V^{h_2})(\phi_1\psi)\\
+&(\dirac^{h_2}+V^{h_0}+V^{h_1})(\phi_2\psi)\\
\end{aligned}
\end{equation*}
for $\psi\in\defD(\dirac^{h})$.
\end{defn}

It is easily verified that the definition is independent of the partition of unity $1=\phi_0+\phi_1+\phi_2$.
\begin{thm}\label{thm:sa}
The Dirac operator $\dirac^h$ is self-adjoint.
\end{thm}
For the proof, we need some lemmas.
\begin{lemma}
The Dirac operator $\dirac:(L_2)^2\to (L_2)^2$ without any magnetic field is a self-adjoint operator with the Sobolev space $(H^1)^2$ as domain.
\end{lemma}

\begin{prf}
See~\cite{th}.\qed
\end{prf}

\begin{lemma}
The Dirac operator $\dirac^{h_0}$ corresponding to the magnetic field $B_0$ is self-adjoint with domain $(H^1)^2$.
\end{lemma}

\begin{prf}
The operator $\dirac^{h_0}$ can be written as $\dirac^{h_0}=\dirac+V^{h_0}$ and the multiplication operator $V^{h_0}$ is relatively bounded with respect to $\dirac$ with relative bound zero, so the lemma follows from the Kato-Rellich theorem.\qed
\end{prf}

\begin{lemma}\label{lem:tozero}
Let $T$ be a bounded operator from $(L_2)^2$ to $(H^1)^2$ and let $V$ be a function, $V(z)\to0$ as $|z|\to\infty$. Then the composition $VT$ is compact in $(L_2)^2$.
\end{lemma}

\begin{prf}
For $n=1,2,\ldots$ we write $V$ as $V=V_n+\tilde{V}_n$, where
\begin{equation*}
V_n(z)=\begin{cases}
V(z) & |V(z)|>\frac{1}{n}\\
0 & |V(z)|\leq \frac{1}{n}.
\end{cases}
\end{equation*}
The functions $V_n$ all have compact support, so the operators $V_nT$ are compact. But $||V_nT-VT||\leq \frac{1}{n}||T||$ for all $n=1,2,\ldots$, so $VT$ is also compact.
\qed
\end{prf}

\begin{remark}
Lemma~\ref{lem:tozero} is also true for $2\times2$ matrix valued functions $V$ where all components tend to zero at infinity. It also holds if $T$ is bounded from $L_2$ to $H^1$.
\end{remark}

\begin{lemma}\label{lem:compact}
Let $0\neq s\in\mathbb{R}$ and let $\phi\in C^\infty_0(\mathbb{R}^2)$ with zero in its support. Then the operator $\phi R$ is compact, where $R=(\dirac^{\tau}+is)^{-1}$ and $\dirac^{\tau}$ is any self-adjoint extension of the Dirac operator corresponding to one AB solenoid (which is assumed to be located at the origin).
\end{lemma}

\begin{prf}
First, $\phi R$ is compact if and only if $\phi R (\phi R)^*=\phi RR^*\phi$ is compact. To show that $\phi RR^*\phi$ is compact, it is sufficient to show that $\phi RR^*$ is compact.

The operator $RR^*$ is equal to $\left((\dirac^\tau)^2+s^2\right)^{-1}$. Note that $(\dirac^\tau)^2$ is a self-adjoint Pauli operator corresponding to the same magnetic field (see Section~\ref{sec:pauli} for a discussion of the Pauli operators that are the square of some Dirac operator). If we denote by $\pauli$ any other self-adjoint Pauli operator corresponding to this magnetic field, then by the Krein resolvent formula (see~\cite{akgl}) the resolvents of $(\dirac^\tau)^2$ and $\pauli$ differ by a finite rank operator. Thus, it is enough to show that $\phi(\pauli+s^2)^{-1}$ is compact for a convenient choice of self-adjoint Pauli extension $\pauli$. Let us choose $\pauli$ to be the Friedrichs extension. The functions in the domain of this extension $\pauli$ vanish at the origin so
\begin{equation*}
\pauli=
\begin{pmatrix}
H & 0\\
0 & H
\end{pmatrix},
\end{equation*}
where $H$ is the Friedrichs extension of the Schrödinger operator corresponding to the same magnetic field (see~\cite{gest} for a discussion of this). Hence it is enough to show that $\phi(H+s^2)^{-1}$ is compact.

Let $H_0=-\Delta$ be the Schrödinger operator corresponding to no magnetic field. Then, by the diamagnetic inequality (see~\cite{meouro}) it follows that
$|\phi(H+s^2)^{-1}u|\leq \phi(H_0+s^2)^{-1}|u|$ (pointwise) for all $u\in L_2$. This inequality implies that $\phi(H+s^2)^{-1}$ is compact if $\phi(H_0+s^2)^{-1}$ is compact (see~\cite{dofr,pi}).

The compactness of $\phi(H_0+s^2)^{-1}$ follows from Lemma~\ref{lem:tozero} since $(H_0+s^2)^{-1}$ is bounded from $L_2$ to $H^1$.\qed
\end{prf}

\begin{lemma}\label{lemma:symmetric}
The operator $\dirac^h$ is symmetric.
\end{lemma}

\begin{prf}
Let $\psi$ and $\tilde{\psi}$ belong to $\defD(\dirac^h)$. Then
\begin{equation*}
\langle \dirac^h\psi,\tilde\psi\rangle = \int (\dirac^h \psi)\cdot \bar{\tilde\psi}d\lambda(z) = \int ((\phi_0+\phi_1+\phi_2)\dirac^h \psi)\cdot \bar{\tilde\psi}d\lambda(z).
\end{equation*}
The symmetry now follows from integration by parts, noticing that (for example)
\begin{equation*}
\int \phi_0 (\dirac^{h_0}+V^{h_1}+V^{h_2})(\phi_0\psi)\cdot\bar{\tilde\psi}d\lambda(z)=\int \phi_0\psi\cdot\overline{(\dirac^{h_0}+V^{h_1}+V^{h_2})(\phi_0\tilde\psi)}d\lambda(z)
\end{equation*}
and
\begin{equation*}
\begin{aligned}
\int \phi_1 (\dirac^{h_0}+V^{h_1}+V^{h_2})(\phi_0\psi)\cdot\bar{\tilde\psi}d\lambda(z)&=\int_{E_{01}}(\dirac^{h_0}+V^{h_1}+V^{h_2})(\phi_0\psi)\cdot\overline{\phi_1\tilde\psi}d\lambda(z)\\
&=\int \phi_0\psi\cdot\overline{(\dirac^{h_1}+V^{h_0}+V^{h_2})(\phi_1\tilde\psi)}d\lambda(z)
\end{aligned}
\end{equation*}
and similar for the other terms. Adding all terms we see that $\dirac^h$ is symmetric.
\qed
\end{prf}

In the following lemma we look at our operator as acting from its domain $\defD(\dirac^h)$ considered as a Hilbert space equipped with graph norm
\begin{equation*}
\begin{aligned}
||\psi||_{\dirac^h}^2&=||(\dirac^{h_0}+V^{h_1}+V^{h_2})(\phi_0 \psi)||^2 + ||(\dirac^{h_1}+V^{h_0}+V^{h_2})(\phi_1 \psi)||^2\\
&+||(\dirac^{h_2}+V^{h_0}+V^{h_1})(\phi_2 \psi)||^2 + ||\psi||^2,
\end{aligned}
\end{equation*}
to $(L_2)^2$.
\begin{lemma}
\label{lemma:fredholm}
Let $0\neq s\in\mathbb{R}$ be fixed. The operator $\dirac^h+is:(\defD(\dirac^h),||\cdot||_{\dirac^h})\to (L_2(\mathbb{R}^2))^2$ is a bounded Fredholm operator with index zero.
\end{lemma}

\begin{prf}
First, it is clear that $\dirac^h+is$ is bounded from the domain space with graph norm. To show that $\dirac^h+is$ is a Fredholm operator, it is enough to find a left and a right parametrix~(see~\cite{agr}). We start by finding a right parametrix. Let $R_j$ denote the resolvent $R_j=(\dirac^{h_j}+is)^{-1}$, $j=0,1,2$, and define the operator $R:(L_2)^2\to(L_2)^2$ as
\begin{equation*}
Ru=\phi_0 R_0u+\phi_1R_1u+\phi_2R_2u,\quad\text{for}\ u\in (L_2)^2.
\end{equation*}
For $u\in (L_2)^2$ we have $\phi_j\phi_kR_ju\in (H^1)^2$ and being zero in a neighborhood of the singular point(s) if $j\neq k$. Thus
\begin{equation*}
(\dirac^{h_k}+V^{h_j})(\phi_j\phi_kR_ju)=(\dirac^{h_j}+V^{h_k})(\phi_j\phi_kR_ju),\quad j\neq k.
\end{equation*}
From this it follows that, if $u\in(L_2)^2$ and $j,k,l\in\{0,1,2,\ j\neq k\neq l\}$, then
\begin{equation*}
(\dirac^h+is)(\phi_jR_ju) = \phi_j u+((V^{h_k}+V^{h_l})\phi_j+\dirac(\phi_j))R_ju.
\end{equation*}
Thus
\begin{equation*}
(\dirac^h+is)Ru = u+K_R u
\end{equation*}
where $K_R:(L_2)^2\to(L_2)^2$ is the operator
\begin{equation*}
\begin{aligned}
K_R u&=\left((V^{h_1}+V^{h_2})\phi_0+\dirac(\phi_0)\right) R_0u\\
&+\left((V^{h_0}+V^{h_2})\phi_1+\dirac(\phi_1)\right) R_1u\\
&+\left((V^{h_0}+V^{h_1})\phi_2+\dirac(\phi_2)\right) R_2u.\\
\end{aligned}
\end{equation*}
$K_R$ is compact. Indeed, the first term is compact according to Lemma~\ref{lem:tozero} since $R_0$ is bounded from $(L_2)^2$ to $(H^1)^2$ and the matrix-valued function $(V^{h_1}+V^{h_2})\phi_0+\dirac(\phi_0)$ tends to zero at infinity. The other two terms are compact by Lemma~\ref{lem:compact}. Hence $K_R$ is compact, so $R$ is a right parametrix.

The operator $L=R_0 \phi_0+ R_1 \phi_1+R_2 \phi_2$ is a left parametrix. Indeed, a calculation similar to the one above shows that

\begin{equation*}
L(\dirac^h+is)\psi= \psi+K_L\psi
\end{equation*}
where
\begin{equation*}
\begin{aligned}
K_L \psi &= R_0 \left((V^{h_1}+V^{h_2})\phi_0-\dirac(\phi_0)\right)\psi\\
&+ R_1 \left((V^{h_0}+V^{h_2})\phi_1-\dirac(\phi_1)\right)\psi\\
&+R_2 \left((V^{h_0}+V^{h_1})\phi_2-\dirac(\phi_2)\right)\psi.
\end{aligned}
\end{equation*}
The compactness of $K_L^*$ follows in the same way as the compactness of $K_R$. Hence $K_L$ is compact. We see that any of $R$ and $L$ works as a parametrix and hence $\dirac^h+is$ is a Fredholm operator.

To see that $\dirac^h+is$ has index zero, we note that since $\dirac$ with domain $(H^1)^2$ is self-adjoint, it has index zero and $R_s:=(\dirac+is)^{-1}$ is a parametrix for $\dirac+is$. The operator
\begin{equation*}
R-R_s=\phi_0 R_0+\phi_1 R_1+\phi_2 R_2-R_s
\end{equation*}
is compact. To see this, write $R-R_s$ as
\begin{equation*}
R-R_s=(\phi_0-1)R_0+\phi_1R_1+\phi_2R_2+(R_0-R_s).
\end{equation*}
The first term is compact according to Lemma~\ref{lem:tozero}, the second and third according to Lemma~\ref{lem:compact}. For the last term we note that $R_0-R_s=-R_0V^{h_0}R_s$. The compactness of $V^{h_0}R_s$ follows from Lemma~\ref{lem:tozero}. Composition with the bounded operator $R_0$ preserves compactness. Thus $R-R_s$ is compact. It follows that
\begin{equation*}
\ind(R)=\ind(R_s).
\end{equation*}
Since $R$ and $R_s$ are parametrices for $\dirac^h+is$ and $\dirac+is$ respectively, it holds that
\begin{equation*}
\ind(R)+\ind(\dirac^h+is)=\ind(R(\dirac^h+is))=0,
\end{equation*}
and
\begin{equation*}
\ind(R_s)+\ind(\dirac+is)=\ind(R_s(\dirac+is))=0.
\end{equation*}
Hence
\begin{equation*}
\ind(\dirac^h+is)=-\ind(R)=-\ind(R_s)=\ind(\dirac+is)=0.
\end{equation*}
\qed    
\end{prf}

\begin{prf}[of Theorem~\ref{thm:sa}]
We know from Lemma~\ref{lemma:symmetric} that $\dirac^h$ is symmetric, so for $0\neq s\in\mathbb{R}$ we have
\begin{equation*}
||(\dirac^h+is)\psi||^2=||\dirac^h\psi||^2+s^2||\psi||^2
\geq s^2||\psi||^2.
\end{equation*}
It follows that $\dim\ker(\dirac^h+is)=0$. From Lemma~\ref{lemma:fredholm} we have that $\dirac^h+is$ has index zero, so it follows that $\dim\ker((\dirac^h)^*-is)=0$. Choosing $s$ positive and negative respectively gives that the deficiency indices for $\dirac^h$ is $(0,0)$, so $\dirac^h$ is self-adjoint.\qed
\end{prf}

\subsection{Spin flip invariance}
\label{sec:diracspinflip}
Since the particle we are studying moves only in a plane, and the magnetic field is orthogonal to this plane, physically it should be no difference if the sign of the magnetic field is changed. This transformation has to come together with a flip of the spin-up and spin-down components and a normalization of the AB intensities. We say that a self-adjoint extension is spin flip invariant if, after applying these transformations, we end up with a (anti)-unitarily equivalent operator. We will show that there are only two values of $\tau$ for which the Dirac operator $\dirac^{\tau,h}$ is spin flip invariant. Let $\btau=(\tau_1,\ldots,\tau_n)$ and denote the Dirac operator by $\dirac^{\btau,h}$. We will use the linear functionals
\begin{align}
\label{eq:diracboundarygen}
c_{-\alpha_j}^{\pm}(\psi_\pm) &= \lim_{r_j\to 0}r_j^{\alpha_j}\frac{1}{2\pi}\int_{0}^{2\pi}\psi_\pm(r_je^{i\theta_j}) d\theta_j,\ \text{and}\\
\label{eq:diracboundarygentwo}
c_{\alpha_j-1}^{\pm}(\psi_\pm) &= \lim_{r_j\to 0}r_j^{1-\alpha_j}\frac{1}{2\pi}\int_{0}^{2\pi}\psi_\pm(r_je^{i\theta_j}) e^{i\theta_j} d\theta_j.
\end{align}
We define anti-unitarily operator $V:(L_2(\mathbb{R}^2))^2\to(L_2(\mathbb{R}^2))^2$ as the spin-flip operator that takes $(\psip,\psim)^t$ to $(\overline{\psim},\overline{\psip})^t$.
\begin{prop}
\label{prop:spinflip}
The operators $\dirac^{\btau',-h}$ and $\dirac^{\btau,h}$ are anti-unitarily equivalent via the operator $V$ if and only if for all $j=1,\ldots,n$ we have $\tau_j'+\tau_j=\pi$ or $\tau_j'+\tau_j=3\pi$.
\end{prop}

\begin{prf}
Let $\beta_j=1-\alpha_j$ be the normalized AB intensities for the magnetic field $-B$ that corresponds to $\dirac^{\btau',-h}$. A function $\psi$ in the domain of $\dirac^{\btau',-h}$ has the asymptotics
\begin{equation*}
\psi \sim \frac{\mu_j}{2}
\begin{pmatrix}
(1+e^{i\tau_j'})\sig{1-\beta_j}r_j^{\beta_j-1}+O(r_j^{1-\beta_j})\\
(1-e^{i\tau_j'})\sig{\beta_j}r^{-\beta_j}e^{i\theta_j}+O(r_j^{\beta_j})
\end{pmatrix}
\end{equation*}
as $z\to z_j$ for some constant $\mu_j\in\mathbb{C}$. We see that $V\psi$ has the asymptotics
\begin{equation*}
\begin{aligned}
V &\psi\sim 
\frac{\bar{\mu}_j}{2}
\begin{pmatrix}
(1-e^{-i\tau_j'})\sig{1-\alpha_j}r_j^{\alpha_j-1}e^{-i\theta_j}+O(r_j^{1-\alpha_j})\\[5pt]
(1+e^{-i\tau_j'})\sig{\alpha_j}r_j^{-\alpha_j}+O(r_j^{\alpha_j})
\end{pmatrix}
\end{aligned}
\end{equation*}
Applying the functionals~\eqref{eq:diracboundarygen} and~\eqref{eq:diracboundarygentwo} we see that $V\psi$ satisfies
\begin{equation*}
\frac{c_{\alpha_j-1}^+((V\psi)_+)}{c_{-\alpha_j}^-((V\psi)_-)}=i\tan(\tau_j'/2)\frac{\sig{1-\alpha_j}}{\sig{\alpha_j}}
\end{equation*}
and $c_{\alpha_j-1}^-((V\psim)_-)=c_{-\alpha_j}^+((V\psip)_+)=0$, so the requirements that the domain change properly is that 
\begin{equation*}
\tan(\tau_j'/2)=\cot(\tau_j/2),\quad\text{for}\ j=1,\ldots,n.
\end{equation*}
We see that $\tau_j/2$ and $\pi/2-\tau_j'/2$ must differ by a integer multiple of $\pi$. Both $\tau_j$ and $\tau_j'$ belong to the interval $[0,2\pi)$, so the only possibilities are $\tau_j'+\tau_j=\pi$ or $\tau_j'+\tau_j=3\pi$.\qed
\end{prf}

\begin{cor}
The operators $\dirac^{\btau,-h}$ and $\dirac^{\btau,h}$ are anti-unitarily equivalent via the operator $V$ if and only if for all $j=1,\ldots,n$ we have $\tau_j=\pi/2$ or $\tau_j=3\pi/2$.
\end{cor}

\begin{prf}
Take $\tau_j'=\tau_j$ in the previous Proposition.\qed
\end{prf}

If we let $W:(L_2(\mathbb{R}^2))^2\to(L_2(\mathbb{R}^2))^2$ be the operator that takes $(\psip,\psim)^t$ to $(\psim,\psip)^t$ we get some other symmetries if we compose it with the gauge transform that only act on the spin-up component.
\begin{prop}
The operators $\dirac^{\btau',-h}$ and $\dirac^{\btau,h}$ are unitarily equivalent via the operator $W$ composed with a gauge multiplication of $e^{-2i\sum\theta_j}$ of the spin-up component if and only if $|\tau_j'-\tau_j|=\pi$ for all $j=1,\ldots,n$.
\end{prop}

\begin{prf}
The proof goes on as in the proof of Proposition~\ref{prop:spinflip}. This time the requirement on $\tau_j$ and $\tau_j'$ becomes
\begin{equation*}
-\tan(\tau_j'/2)=\cot(\tau_j/2),\quad\text{for}\ j=1,\ldots,n
\end{equation*}
which gives $|\tau_j'-\tau_j|=\pi$ for all $j=1,\ldots,n$.\qed
\end{prf}

\subsection{Zero-modes}
Let us calculate the dimension of the kernel of $\dirac^h$ under the assumption that $\tau_j=\tau$ for all $j=1,\ldots,n$, which means that we assume that we have the same physical conditions of the behavior of the particle close to all solenoids. Denote by $\Phi$ the total flux of $B$ divided by $2\pi$, that is
\begin{equation*}
\Phi=\frac{1}{2\pi}\int B(z)d\lambda(z)=\frac{1}{2\pi}\int B_0(z)d\lambda(z)+\sum_{j=1}^n\alpha_j.
\end{equation*}
As usual, the definition of the total flux is a matter of agreement, due to the arbitrariness in the choice of normalization for AB intensities. The asymptotics of $e^h$ at infinity and at the singular points $\Lambda$ are given by
\begin{equation}
\label{eq:asymp}
e^h\sim
\begin{cases}
|z|^\Phi & |z|\to\infty\\
|z-z_j|^{\alpha_j} & z\to z_j.
\end{cases}
\end{equation}
Remember also that the functions in the domain of $\dirac^h$ must satisfy
\begin{equation}
\label{eq:singasymp}
\frac{c_{\alpha_j-1}^+(\psip)}{c_{-\alpha_j}^-(\psim)}=i\cot(\tau_j/2)\frac{\sig{1-\alpha_j}}{\sig{\alpha_j}},\quad j=1,\ldots,n.
\end{equation}
Let $\{x\}$ denote the lower integer part, that is
\begin{equation*}
\{x\}=
  \begin{cases}
\lfloor x\rfloor & x>1\ \text{and}\ x\not\in\mathbb{N}\\
x-1 & x>1\ \text{and}\ x\in\mathbb{N}\\
0 & \text{otherwise}.
  \end{cases}
\end{equation*}

\begin{thm}
\label{thm:ACDirac}
If $\tau_j=\tau$, $j=1,\ldots,n$ then the dimension of the kernel of $\dirac^h$ is given by
\begin{equation*}
\dim\ker \dirac^h =
\begin{cases}
\{|n-\Phi|\}&\text{if}\ \tau=0,\\
\{|\Phi|\}&\text{if}\ \tau=\pi,\\
0&\text{otherwise}.
\end{cases}
\end{equation*}
\end{thm}
The proof follows the same idea as the original proof by Aharonov-Casher with the same changes as in~\cite{per} and using the fact that the spin-up and spin-down components are coupled if $\tau\not\in\{0,\pi\}$.
\begin{prf}
We start by calculating the zero-modes as if the spin-up and spin-down components were not coupled; so these components are studied separately.

Let us start with the spin-up component, that is, we consider the solutions to $\qplus\psip=0$. This is equivalent to $\dzbar(e^{-h}\psip)=0$, and thus the function $f_+=e^{-h}\psip$ must be analytic in $\mathbb{C}\setminus\Lambda$. The behavior of $f_+$ at the singular points $\Lambda$ is different for different values of the parameter $\tau$, but a pole of order at most $\{-\Phi\}-1$ at infinity is allowed independently of the value of $\tau$.

{\bf Case I}, $\tau=\pi$: For $\psip$ to belong to $L_2$, as we see from~\eqref{eq:asymp}, the function $f_+$ is not allowed to have any poles at the singular points $\Lambda$. Thus, if $\tau=\pi$ then $f_+$ may be a polynomial of order at most $\{-\Phi\}-1$. There are as many as $\{-\Phi\}$ many linearly independent such polynomials.

{\bf Case II}, $\tau\neq\pi$: From~\eqref{eq:singasymp} we see that a pole of order at most one is allowed at each $z_j\in\Lambda$. The calculation in~\cite{per} then yields that the dimension is $\{n-\Phi\}$.

Let us now turn to the spin-down component. We look for solutions to the equation $\qminus\psim=0$, which is equivalent to finding solutions to $\dz(e^h\psim)=0$. If we now let $f_-=e^h\psim$, then $f_-$ must be anti-analytic in $\mathbb{C}\setminus\Lambda$, and from the asymptotics~\eqref{eq:asymp} we see that $f_-$ may have a polynomial part of degree at most $\{\Phi\}-1$ independent of the value of the parameter $\tau$. Again we get two different cases for the behavior of the functions at the singular points $\Lambda$.

{\bf Case I}, $\tau=0$: In this case we see from~\eqref{eq:singasymp} that no singular parts for $\psim$ are allowed at $\Lambda$, and hence $f_-$ must have a zero of order at least 1 at each point in $\Lambda$. That is we have a polynomial in $\bar{z}$ of degree $\{\Phi\}-1$ with $n$ predicted zeroes. There are $\{\Phi-n\}$ linearly independent polynomials of this type.

{\bf Case II}, $\tau\neq 0$: Now $f_-$ must be a polynomial in $\bar{z}$ of degree at most $\{\Phi\}-1$, but without any forced zeroes. Thus the dimension of the kernel is $\{\Phi\}$. 

Since the spin-up and spin-down components are not coupled in the cases $\tau=0$ and $\tau=\pi$ the calculations above yield
\begin{equation*}
\dim\ker \dirac^h =
\begin{cases}
\{|n-\Phi|\}&\text{if}\ \tau=0,\\
\{|\Phi|\}&\text{if}\ \tau=\pi.\\
\end{cases}
\end{equation*}
Let us now assume that $\tau\not\in\{0,\pi\}$. We should evaluate how the spin-up zero-modes match the spin-down zero-modes to satisfy the conditions at the singularities. First we note that to be able to have zero-modes both $\{n-\Phi\}$ and $\{\Phi\}$ must be positive. From the calculations in the last two paragraphs of the proof of Theorem~3.3 in~\cite{per} it follows that $f_+$ must be of the form
\begin{equation*}
f_+(z)=\sum_{j=1}^n \frac{\eta_j}{z-z_j}
\end{equation*}
where $\eta_j\in\mathbb{C}$ satisfy
\begin{equation}
\label{eq:cjkrav}
\sum_{j=1}^n \eta_j z_j^k=0,\quad \text{for}\ k=0,1,\ldots,n-\{n-\Phi\}-1
\end{equation}
and $f_-(z)$ must be a polynomial in $\bar{z}$ of degree at most $\{\Phi\}-1$. Actually, we will show that even if the degree of the polynomial $f_-$ is $\{\Phi\}$ or in some cases $\{\Phi\}+1$, all coefficients of the polynomial must be zero. Let us define the natural number $m$ as $m=n-\{n-\Phi\}-1$ and note that $m=\lfloor \Phi\rfloor$. Let
\begin{equation}
\label{eq:polynom}
f_-(z)=\sum_{j=0}^{m}a_j \bar{z}^j.
\end{equation}
From the asymptotics~\eqref{eq:singasymp} we see that
\begin{equation*}
\frac{\eta_j}{f_-(z_j)}=e^{-2h_0(z_j)}\prod_{l\neq j}\left(|z_j-z_l|^{-2\alpha_l}\right)i\cot(\tau/2)\frac{\sig{1-\alpha_j}}{\sig{\alpha_j}},\quad j=1,\ldots,n.
\end{equation*}
From the requirements~\eqref{eq:cjkrav} of the coefficients $\eta_j$ we get
\begin{equation}
\label{eq:nollekv}
0=\sum_{j=1}^n \eta_j z_j^k=i\cot(\tau/2)\sum_{j=1}^n b_j f_-(z_j) z_j^k,\quad k=0,1,\ldots,m,
\end{equation}
where
\begin{equation*}
b_j=e^{-2h_0(z_j)}\prod_{l\neq j}\left(|z_j-z_l|^{-2\alpha_l}\right)\frac{\sig{1-\alpha_j}}{\sig{\alpha_j}}>0.
\end{equation*}
We introduce the vector $\mathbf{a}=(a_0,\ldots,a_{m})^t$ where $a_j$ are the coefficients in~\eqref{eq:polynom}. Let us also introduce the matrix $V$ as
\begin{equation*}
V=
\begin{pmatrix}
1 & 1 & \cdots & 1\\
z_1 & z_2 & \cdots & z_n\\
\vdots & \vdots & \ddots & \vdots\\
z_1^m & z_2^m & \cdots & z_n^m\\
\end{pmatrix},
\end{equation*}
and the diagonal matrix $B$ having the positive number $b_j$ at the $j$th diagonal position. Then~\eqref{eq:nollekv} can be written as
\begin{equation*}
i\cot(\tau/2)VBV^*\mathbf{a}=0.
\end{equation*}
The matrix $VBV^*$ is clearly Hermitian and since $B$ is positive, we can write $VBV^*$ as $(V\sqrt{B})(V\sqrt{B})^*$. Hence the null space of $VBV^*$ is the same as that of the matrix $(V\sqrt{B})^*=\sqrt{B}V^*$. Since $V^*$ is (a part of) a Vandermonde matrix it has full rank, so the dimension of the null space of $\sqrt{B}V^*$ is zero. Hence the polynomial $f_-$, and thus also $\psim$, must be zero. Since the spin-up and spin-down components are coupled, it follows that $\psip$ is also zero. Consequently, $\dim\ker\dirac^h=0$, and the proof is complete.\qed
\end{prf}

\section{The Pauli operator}
\label{sec:pauli}

In this section we will study the Pauli operator corresponding to the magnetic field~\eqref{eq:magn}, obtained as the square of a self-adjoint Dirac operator from the previous section. 

\begin{defn}
\label{def:pauli}
We define the Pauli operator $\pauli^h$ as $(\dirac^h)^2$ where $\dirac^h$ is a self-adjoint Dirac operator defined in Definition~\ref{def:dirac}. This means that
\begin{equation*}
\begin{aligned}
\defD(\pauli^h) &= \{\psi\in (L_2)^2:\ \diraca^h\psi\in\defD(\dirac^h)\},\\
\pauli^h\psi &= (\diraca^h)^2\psi,\quad\text{for}\ \psi\in\defD(\pauli^h).
\end{aligned}
\end{equation*}
\end{defn}
Let us again introduce the boundary value linear functionals acting on $\defD(\pauli^h)$, but this time for all singular points $\Lambda$. For $j=1,\ldots,n$, let
\begin{align*}
c_{-\alpha_j}^{\pm}(\psi_\pm) &= \lim_{r_j\to 0}r_j^{\alpha_j}\frac{1}{2\pi}\int_{0}^{2\pi}\psi_\pm(r_je^{i\theta_j}) d\theta_j,\\
c_{\alpha_j}^{\pm}(\psi_\pm) &= \lim_{r_j\to 0}r_j^{-\alpha_j}\left(\frac{1}{2\pi}\int_{0}^{2\pi}\psi_\pm(r_je^{i\theta_j}) d\theta_j-r_j^{-\alpha_j}c_{-\alpha_j}^{\pm}(\psi_\pm)\right),\\
c_{\alpha_j-1}^{\pm}(\psi_\pm) &= \lim_{r_j\to 0}r_j^{1-\alpha_j}\frac{1}{2\pi}\int_{0}^{2\pi}\psi_\pm(r_je^{i\theta_j}) e^{i\theta_j} d\theta_j,\ \text{and}\\
c_{1-\alpha_j}^{\pm}(\psi_\pm) &= \lim_{r_j\to 0}r_j^{\alpha_j-1}\left(\frac{1}{2\pi}\int_{0}^{2\pi}\psi_\pm(r_je^{i\theta_j}) e^{i\theta_j}d\theta_j-r_j^{\alpha_j-1}c_{\alpha_j-1}^{\pm}(\psi_\pm)\right).
\end{align*}
Since there are more self-adjoint Pauli extensions than Dirac extensions corresponding to our singular magnetic field, it is clear that not all Pauli operators can be obtained as the square of a self-adjoint Dirac operator.
\begin{prop}
For an arbitrary self-adjoint Pauli extension $\pauli$, it is the square of some self-adjoint Dirac extension $\dirac^h$ if and only if the following equations are satisfied for all $\psi\in\defD(\pauli)$
\begin{align}
\label{eq:pcondf}
\frac{c_{\alpha_j-1}^+(\psip)}{c_{-\alpha_j}^-(\psim)}&=i\cot(\tau_j/2)\frac{\sig{1-\alpha_j}}{\sig{\alpha_j}},\\
\label{eq:pauliextra}
\frac{c_{\alpha_j}^-(\psim)}{c_{1-\alpha_j}^+(\psip)}&=i\cot(\tau_j/2)\frac{\sig{-\alpha_j}}{\sig{\alpha_j-1}},\\
c_{-\alpha_j}^+(\psip) &=0,\ \text{and}\\
\label{eq:pcondl}
c_{\alpha_j-1}^-(\psim) &=0.
\end{align}
\end{prop}

\begin{prf}
Let $\dirac^h$ be a given self-adjoint Dirac extension, and let $\psi$ belong to $\defD(\dirac^h)$. Then for some constants $\mu_j$ we have
\begin{equation*}
\psi\sim
\frac{\mu_j}{2}
\begin{pmatrix}
(1+e^{i\tau_j})\sig{1-\alpha_j}r_j^{\alpha_j-1}e^{-i\theta_j}+O(r_j^{\alpha_j})\\[5pt]
(1-e^{i\tau_j})\sig{\alpha_j}r_j^{-\alpha_j}+O(r_j^{1-\alpha_j})
\end{pmatrix},\quad \text{for}\ j=1,\ldots,n
\end{equation*}
We want to find the next term in the asymptotical expansion for $\psi$ such that
\begin{equation*}
\dirac^h(\psi)\sim
\frac{\nu_j}{2}
\begin{pmatrix}
(1+e^{i\tau_j})\sig{1-\alpha_j}r_j^{\alpha_j-1}e^{-i\theta_j}\\[5pt]
(1-e^{i\tau_j})\sig{\alpha_j}r_j^{-\alpha_j}
\end{pmatrix},\quad j=1,\ldots,n
\end{equation*}
for some constants $\nu_j$. This means that $\psi$ must have the asymptotics
\begin{small}
\begin{equation*}
\psi\sim 
\begin{pmatrix}
\frac{\mu_j}{2}(1+e^{i\tau_j})\sig{1-\alpha_j}r_j^{\alpha_j-1}e^{-i\theta_j}-\frac{i\nu_j}{2}(1-e^{i\tau_j})\sig{\alpha_j-1}r_j^{1-\alpha_j}e^{-i\theta_j}+O(r_j^{2-\alpha_j})\\[5pt]
\frac{\mu_j}{2}(1-e^{i\tau_j})\sig{\alpha_j}r_j^{-\alpha_j}-\frac{i\nu_j}{2}(1+e^{i\tau_j})\sig{-\alpha_j}r^\alpha_j+O(r_j^{1+\alpha_j})
\end{pmatrix}
\end{equation*}
\end{small}
as $z$ tends to $z_j$, for $j=1,\ldots,n$. From this it follows that
\begin{align}
\label{eq:coefff}
c_{\alpha_j-1}^+(\psip) &=\frac{\mu_j}{2}(1+e^{i\tau_j})\sig{1-\alpha_j},\\
c_{1-\alpha_j}^+(\psip) &=-\frac{i\nu_j}{2}(1-e^{i\tau_j})\sig{\alpha_j-1},\\
c_{-\alpha_j}^-(\psim) &=\frac{\mu_j}{2}(1-e^{i\tau_j})\sig{\alpha_j},\ \text{and}\\
c_{\alpha_j}^-(\psim) &= -\frac{i\nu_j}{2}(1+e^{i\tau_j})\sig{-\alpha_j}.
\end{align}
Moreover
\begin{align}
c_{-\alpha_j}^+(\psip)&=0,\ \text{and}\\
c_{\alpha_j-1}^-(\psim)&=0
\end{align}
since no such singular functions belong to $\defD(\dirac^h)$. The coefficients
\begin{align}
c_{\alpha_j}^+(\psip)&=\text{arbitrary},\ \text{and}\\
\label{eq:coeffl}
c_{1-\alpha_j}^-(\psim)&=\text{arbitrary}
\end{align}
since such terms disappear (near the singular point) when applying $\dirac^h$.\qed
\end{prf}

\begin{remark}
The definition of $\pauli^h$ can be written as
\begin{equation*}
\begin{aligned}
\defD(\pauli^h) &= \Big\{\psi\in(L_2(\mathbb{R}^2))^2:\
\diraca^2 \psi\in(L_2(\mathbb{R}^2))^2;\\
&\quad\quad \text{equations~\eqref{eq:pcondf}--\eqref{eq:pcondl} hold for all $\psi$}\Big\}.
\end{aligned}
\end{equation*}
We see also that $\defD(\pauli^h)$ is exactly the subset of $\defD(\dirac^h)$ for which also the condition~\eqref{eq:pauliextra} holds.\qed
\end{remark}

\subsection{Spin-flip invariance and Zero-modes}

\begin{prop}
The only self-adjoint Pauli extensions $\pauli^{\btau,h}=(\dirac^{\btau,h})^2$ that are spin-flip invariant under the transform V are these where for all $j=1,\ldots,n$ we have $\tau_j=\pi/2$ or $\tau_j=3\pi/2$.
\end{prop}

\begin{prf}
The proof is the same as for the Dirac operators, see Proposition~\ref{prop:spinflip}.\qed
\end{prf}

\begin{thm}
\label{thm:ACPauli}
If $\tau_j=\tau$, $j=1,\ldots,n$ then the dimension of the kernel of $\pauli^h$ is given by
\begin{equation*}
\dim\ker \pauli^h =
\begin{cases}
\{|n-\Phi|\}&\text{if}\ \tau=0,\\
\{|\Phi|\}&\text{if}\ \tau=\pi,\\
0&\text{otherwise}.
\end{cases}
\end{equation*}
\end{thm}

\begin{prf}
This is a direct consequence of Theorem~\ref{thm:ACDirac} since $\ker \pauli^h=\ker \dirac^h$.\qed
\end{prf}

\subsection{Discussion}\label{sec:disc}
Let us compare the different self-adjoint Pauli operators from~\cite{ervo} (which we will denote by $\evpauli$) and~\cite{per} (which we will denote by $\mpauli$) with the ones obtained above as the square of a self-adjoint Dirac operator. It is easier to do this comparison if we have the same AB flux normalization for all operators. Thus, we let all AB intensities $\alpha_j$ belong to the interval $(0,1)$. In the case of the Pauli operator $\evpauli$, where the AB intensities were normalized to $[-1/2,1/2)$, we have to do a gauge transformation if there are intensities $\alpha_j$ belonging to $[-1/2,0)$. This is not a problem, since $\evpauli$ is gauge invariant.

In Table~\ref{tab:pauli} we see a comparison of the boundary conditions of the Pauli operators obtained above that are the square of a Dirac operator and the Maximal and EV Pauli operators (see~\cite{per,ervo}). We see that $\mpauli$ is not the square of a Dirac operator. However, if we let
\begin{equation*}
\tau_j=
\begin{cases}
\pi, & \text{if}\quad 0<\alpha_j<1/2\\ 
0, & \text{if}\quad 1/2\leq\alpha_j<1
\end{cases},\quad j=1,\ldots,n,
\end{equation*}
and $\btau=(\tau_1,\ldots,\tau_n)$, then $\evpauli$ is the square of the self-adjoint Dirac operator corresponding to $\btau$. Note that it is possible to have different physical situations at the singular points $\Lambda$. Indeed, if not all intensities $\alpha_j$ belong to either $(0,1/2)$ or $[1/2,1)$ then this is the case.

\begin{table}[!ht]
\caption{\label{tab:pauli}The boundary value conditions for the squared Dirac operators compared with the ones for the Maximal and EV Pauli operators.}
\begin{center}
\begin{tabular}{l|lll}
 & $\pauli^h=(\dirac^h)^2$ & $\mpauli$ & $\evpauli$ \\
\hline
$\frac{c_{\alpha_j}^+}{c_{-\alpha_j}^+}$      & $\infty$ & $\infty$ & $\infty$ \\
\hline\\[-3mm]
$\frac{c_{1-\alpha_j}^+}{c_{\alpha_j-1}^+}$ & $-\frac{\nu_j}{\mu_j}\frac{\sig{\alpha_j-1}}{\sig{1-\alpha_j}}\tan(\tau_j/2)$ & $0$ & \begin{small}$\begin{cases}\infty, & \text{if}\quad 0<\alpha_j<1/2\\ 0, & \text{if}\quad 1/2\leq\alpha_j<1\end{cases}$\end{small} \\
\\[-3mm]
\hline\\[-3mm]
$\frac{c_{\alpha_j}^-}{c_{-\alpha_j}^-}$    & $\frac{\nu_j}{\mu_j}\frac{\sig{-\alpha_j}}{\sig{\alpha_j}}\cot(\tau_j/2)$ & $0$ &  \begin{small}$\begin{cases}0, & \text{if}\quad 0<\alpha_j<1/2\\ \infty, & \text{if}\quad 1/2\leq\alpha_j<1\end{cases}$\end{small}\\
\\[-3mm]
\hline
$\frac{c_{1-\alpha_j}^-}{c_{\alpha_j-1}^-}$  & $\infty$ & $\infty$ & $\infty$
\\
\hline
\end{tabular}
\end{center}
\end{table}
\begin{remark}
If the AB intensities in~\cite{ervo} would have been normalized to $(0,1)$ instead of $[-1/2,1/2)$, then the operator $\evpauli$ would have become the square of the Dirac operator where $\tau_j=\pi$ for all $j=1,\ldots,n$. If the AB intensities would have been normalized to $(-1,0)$ then $\evpauli$ would have been the square of the Dirac operator where $\tau_j=0$ for all $j=1,\ldots,n$.\qed
\end{remark}
Among the Pauli operators studied in this article, the ones for $\tau=\pi/2$ (which is (anti)-unitarily equivalent to the one for $\tau=3\pi/2$), $\tau=0$ and $\tau=\pi$ seems to be the most interesting ones. For $\tau=\pi/2$ we get a very symmetric domain of the operator, which implies that the operator is spin-flip invariant. Lacking zero-modes, it does not satisfy the original Aharonov-Casher formula, but it can be approximated component-wise according to Table~\ref{tab:pauli} and the result in~\cite{bopu}. See the end of~\cite{per} for a discussion of this.

The Pauli operators corresponding to $\tau=0$ and $\tau=\pi$ have very asymmetric domains. Only one of the components contain singular terms at the points $\Lambda$. This lack of symmetry implies that these extensions are not spin-flip invariant. On the other hand, the Pauli operator corresponding to $\tau=\pi$ does satisfy the original Aharonov-Casher formula and there is no doubt that both of these Pauli operators can be approximated as in~\cite{bopu}, even as Pauli Hamiltonians.

The Maximal Pauli operator studied~\cite{per} is spin-flip invariant and has zero-modes, even more than is present in the original Aharonov-Casher formula. It can be approximated component-wise as in~\cite{bopu}. However, it is not the square of a self-adjoint Dirac operator.

It is still not clear which Pauli extension that describes the physics in the best way.

\section*{Acknowledgments}
I would like to thank my supervisor, Professor Grigori Rozenblum, for assisting me during the work and coming up with the idea of the proof of Lemma~\ref{lem:compact}.

\bibliography{mp}
\bibliographystyle{plain}
\end{document}